\documentclass[letterpaper, 10 pt, conference]{ieeeconf}  

\IEEEoverridecommandlockouts                              
\usepackage{algorithm}
\usepackage{algorithmic}
\usepackage{xcolor}
\usepackage{lipsum}
\usepackage{hyperref}

\usepackage{graphics} 
\usepackage{epsfig} 
\usepackage{mathptmx} 
\usepackage{times} 
\usepackage{amsmath} 
\usepackage{amssymb}  

\usepackage{tikz}
\usetikzlibrary{plotmarks}
\usepackage{caption}
\usepackage{subcaption}
\usepackage{pgfplots}
\usepgfplotslibrary{fillbetween}
\pgfplotsset{compat=1.18}
\usepackage{comment}
\usepackage{xcolor}
\usepackage{pifont}
\usepackage{pgfplots}
\usepackage{ulem}
\usepackage{color}
\usepackage{float}
\usepackage{stfloats}
\usepackage{cite}

\usepackage{siunitx}
\usepackage{pgfplotstable}

    \pgfplotsset{
        compat=1.16,
        layers/MyLayers/.define layer set={
            waybackgroundlayer, boundingboxlayer,
            axis background,
            pre main,axis grid,axis ticks,axis lines,axis tick labels,
            axis descriptions,axis,main,foreground,
        }{/pgfplots/layers/axis on top},
        set layers=MyLayers,
    }

\title{\LARGE \bf
Mean Field Control of Thermostatically Controlled Loads \\as Piecewise Deterministic Markov Processes
}

\author{ \parbox{2 in}{\centering Thomas Le Corre\\
         Inria, Paris, France \\ DI ENS, ENS, PSL University \\ Paris, France\\
         {\tt\small thomas.le-corre@inria.fr}}
         \hspace*{ 0.3 in}
         \parbox{2 in}{ \centering Adrien Séguret\\
         OSIRIS department, EDF Lab\\
         Palaiseau, France\\
         {\tt\small adrien.seguret@edf.fr}}
         \hspace*{ 0.3 in}
         \parbox{2 in}{ \centering Ana Bušić\\
        Inria, Paris, France \\ DI ENS, ENS, PSL University \\ Paris, France\\
         {\tt\small ana.busic@inria.fr}}
}

\DeclareMathAlphabet{\mathcal}{OMS}{cmsy}{m}{n}

\begin{document}

\maketitle
\thispagestyle{empty}
\pagestyle{empty}

\renewcommand{\thefootnote}{} 
\footnotetext{This work was carried out in the framework of the AI-NRGY project, funded by France 2030 (Grant No: ANR-22-PETA-0004).}
\renewcommand{\thefootnote}{\arabic{footnote}}
\begin{abstract}
This paper presents a mean-field control approach for \textit{Piecewise Deterministic Markov Processes} (PDMPs), specifically designed for controlling a large number of agents. By modeling the interactions of a large number of agents through an aggregate cost function, the proposed method mitigates the high dimensionality of the problem by focusing on a representative agent. 
{ The contribution of this work is the application of a PDMP-based mean-field control framework to the coordination of a large population of Thermostatically Controlled Loads (TCLs). Adapting this framework to TCLs requires incorporating a quality-of-service constraint ensuring that each agent’s temperature remains within a specified comfort range. To achieve this, an additional jump intensity is introduced so that agents are very likely to switch between heating and cooling modes when they reach the boundaries of their temperature range.}
This extension to TCLs is demonstrated through \textit{Water Heaters} (WHs) control, with a decentralized algorithm based on a dual formulation and stochastic gradient descent. The numerical results obtained illustrate this approach on two examples (signal tracking and taking into account energy price).
\end{abstract}

\section{INTRODUCTION}

Mean-field control problems have gained significant attention due to their numerous applications in economics, energy management, and network engineering. When a large number of agents interact through an aggregate cost function, direct resolution becomes impractical due to the dimensionality of the problem. The mean-field approximation mitigates this complexity by modeling a representative agent evolving under the influence of the population distribution. This approach has led to agent control problems, where the objective is to coordinate them with a common goal like demand response for TCLs.

Piecewise deterministic Markov processes (PDMPs) provide a relevant framework for modeling such dynamics, particularly when dealing with systems involving sudden transitions and deterministic evolution between these transitions \cite{davis1993markov}, and applications to TCLs have already been proposed \cite{BusicMeyn2019}. Recently, a decentralized approach based on a Lagrangian formulation was proposed to solve these problems \cite{seguret2024decentralized}. This approach relies on a stochastic Uzawa-type algorithm to efficiently handle the non-convex setting.

This research find applications in the need for demand-side management in power systems to ensure a balance between energy production and consumption. With the growing integration of intermittent renewable sources such as wind and solar, managing electricity demand has become essential to counterbalance the variability of power generation.  The literature on demand-side management and distributed control has investigated various strategies \cite{hochberg2006demand,garabe2022ensemble}. Previous studies have successfully addressed this challenge by controlling the aggregate consumption of flexible consumers \cite{BusicMeyn2019}. Additionally, mean-field approaches have been explored to optimize the coordinated charging of large fleets of electric vehicles (EVs), formulating the problem as an optimal control of partial differential equations modeling PDMPs \cite{seguret2024decentralized}. Mean-field models for TCLs have long been considered \cite{malhame1985,malhame1988,laurent1994}, with individual TCL modeled by Markov Jump Process \cite{malhame1990}. Various mean field approaches have proved to be successful, whether by modeling the TCLs as price-responsive rational agents\cite{depaola2018}, for controlling their mean temperature \cite{li2016,kizilkale2014} or to lower power peaks \cite{grammatico2015}.

\textbf{Contributions:}
{
In this paper, we adapt a PDMP-based mean-field control framework to the coordination of a large population of Thermostatically Controlled Loads (TCLs).
\begin{itemize}
    \item For safety or quality-of-service reasons, the temperature of a TCL must remain within a desired range, so we introduced an additional jump intensity to switch when reaching the boundary of this range.
    \item We illustrate this approach by a case study of water heater control to show that it can be used for various goals, such as signal tracking or reducing total electricity cost.
\end{itemize}
}

\textbf{Structure:}
The paper is structured as follows. Section \ref{s:Pre} introduces the mathematical tools needed, and in particular the definition of PDMPs. Section \ref{s:Prob} recalls the mean field control problem and its dual formulation, the main assumptions, and the main theoretical results. Finally, Section \ref{s:Alg} introduces the Uzawa algorithm and Section \ref{s:WH} illustrates the performance of our method in a concrete application, where a central planner aims to control a large population of water heaters.

\section{PRELIMINARIES}
\label{s:Pre}
\textbf{Notations:} To facilitate readability, we introduce the main notations used throughout the paper.

\begin{itemize}
    \item Let $\mathcal{X}$ be the state space defined as $\mathcal{X} := I \times \mathbb{R}$, where $I$ is a finite set with cardinality $d \in \mathbb{N}^*$.
    \item The time horizon is given by the interval $[0,T]$ with $T > 0$.
    \item The space $C^1(\mathcal{X})$ (resp. $C^1([0,T] \times \mathcal{X})$) represents the set of continuously differentiable real-valued functions defined on $\mathcal{X}$ (resp. $[0,T] \times \mathcal{X}$).
    \item The space $C_0([0,T] \times \mathcal{X}, \mathbb{R}^d_+)$ denotes the set of continuous functions from $[0,T] \times \mathcal{X}$ to $\mathbb{R}^d_+$, equipped with the supremum norm:
    \begin{equation*}
        \|f\|_{\infty} := \sup_{(i,t,j,\theta) \in I \times [0,T] \times \mathcal{X}} |f_i(t,j,\theta)|.
    \end{equation*}
    \item We define the following function space:
    
    \begin{align*}
        A &:= \{ \alpha \in C_0([0,T] \times \mathcal{X}, \mathbb{R}^d_+) \mid \forall i \in I, \alpha_i(\cdot, i, \cdot) = 0 \}.
    \end{align*}
    \item The space of Borel probability measures on $\mathcal{X}$ is denoted by $\mathcal{P}(\mathcal{X})$.
    \item The space $C_0([0,T], \mathcal{P}(\mathcal{X}))$ consists of continuous functions from $[0,T]$ to $\mathcal{P}(\mathcal{X})$, endowed with the metric:
    \begin{equation*}
        W_{\| \cdot \|_{\infty}}(m_1, m_2) := \sup_{t \in [0,T]} W(m_1(t), m_2(t)).
    \end{equation*}
\end{itemize}

\textbf{Definition of PDMP:}
A Piecewise Deterministic Markov Process (PDMP) \cite{davis1993markov} $X_t = (I_t, \theta_t) \in \mathcal{X}$ is defined by:
\begin{itemize}
    \item A deterministic flow between jumps: the continuous component $\theta_t$ { representing the temperature} evolves according to the ordinary differential equation:
    \begin{equation*}
        \frac{d}{dt} \theta_t = b(t,I_t, \theta_t),
    \end{equation*}
    where $b: I \times \mathbb{R} \to \mathbb{R}$ describes the speed of $\theta_t$, { i.e. the dynamics of the temperature}. For any \((\tau,t,j,\theta) \in [0,T] \times [\tau,T] \times \mathcal{X}\), we define the flow \(\phi\) as the unique solution of this ordinary differential equation $\partial_t \phi_{\tau,\theta}(j,t) = b\big(j,\phi_{\tau,\theta}(j,t)\big)$ with the initial condition $\phi_{\tau,\theta}(j,\tau) = \theta$. 
    \item A stochastic jump mechanism: the discrete component $I_t$ {  representing the heating mode}, switches between states in $I$ at jump times $\{T_k\}_{k \geq 0}$, which are determined by a state-dependent intensity function $\alpha: [0,T] \times \mathcal{X} \to \mathbb{R}_+$:
    \begin{equation*}
        P(T_{k+1} - T_k > s | X_{T_k} = (Y,\theta)) = e^{- \int_{0}^{s} \sum_{j \neq Y} \alpha_j(t, Y,\theta) dt}.
    \end{equation*}
    \item At each jump time $T_k$, the process transitions from state $(I_{T_k^-}, \theta_{T_k^-})$ to a new state $(I_{T_k}, \theta_{T_k})$:
    \begin{align*}
        I_{T_k} &= i \quad \text{with probability } \frac{\alpha_i(T_k, I_{T_k^-}, \theta_{T_k^-})}{\sum_{i' \neq I_{T_k^-}} \alpha_{i'}(T_k, I_{T_k^-}, \theta_{T_k^-})}, \\
        \theta_{T_k} &= \theta_{T_k^-}.
    \end{align*}
    { In practice, the sequence $(T_0,T_1,\ldots, T_k,\ldots)$ represents the ON / OFF switching times of the TCL}
\end{itemize}

The PDMP thus alternates between deterministic evolution phases and stochastic jumps, making it a suitable model for systems with sudden transitions and deterministic dynamics between jumps. {  Modeling the TCLs as PDMPs also results in a randomized control $\alpha$, thereby desynchronizing the behavior of the TCLs.}

\section{PROBLEM FORMULATION}
\label{s:Prob}
\subsection{Context}
We consider a population of $N$ independent and identically distributed agents, each following a Piecewise Deterministic Markov Process (PDMP) of jump intensity $\hat{\alpha}+\alpha $. Here, $\hat{\alpha}$ is a function in $A$, and the addition of this function forces jumps with high probability when the agent is in certain states. For example, in the water heater application, we want to maintain the temperature between two bounds: $\forall t, \theta\in [\theta_{min},\theta_{max}]$. With $\hat{\alpha}$ large when $\theta\geq\theta_{max}$ and $i=1$ or when $\theta\leq\theta_{min}$ and $i=0$, we force the PDMP to remain within a range of desirable values. The second term $\alpha\in A$ is our control parameter.

The state of an agent at time $t$ is given by $X_t^{\alpha+\hat{\alpha}} = (I_t^{\alpha+\hat{\alpha}},\theta_t^{\alpha+\hat{\alpha}})$, where:
\begin{itemize}
    \item $I_t^{\alpha+\hat{\alpha}}\in I$ represents a discrete mode (e.g., heating mode for a water heater),
    \item $\theta_t^{\alpha+\hat{\alpha}}\in \mathbb{R}$ represents a continuous state variable (e.g., temperature of a water heater).
\end{itemize}

The goal is to minimize the expected cost function $J_N$ over all admissible controls $\alpha \in A$:
\begin{equation}\label{def_jn}
\begin{aligned}
    J_N&(\alpha) :=\\
    \mathbb{E}& \left[ \int_0^T  f \left( t, \frac{1}{N} \sum_{n=1}^{N} p(t, X_t^{n,{\alpha+\hat{\alpha}}}) \right) dt + \frac{1}{N} \sum_{n=1}^{N} G(\alpha, X^{n,\alpha+\hat{\alpha}}) \right]
\end{aligned}
\end{equation}
where:
\begin{itemize}
    \item $p(t, X_t^{n,\alpha+\hat{\alpha}})$ represents an individual consumption,
    \item $f(t, \cdot)$ models a coupling cost depending on the aggregate consumption of the agents, 
    \item $G(\alpha, X^{n,\alpha+\hat{\alpha}})$ is the individual cost of an agent, given by:
    \begin{equation}
    \label{e:defG}
        G(\alpha, x) := \int_0^T \left[ c(t, x_t) + \sum_{j \in I} L(\alpha_j (t, x_t)) \right] dt + g(x_T).
    \end{equation}
\end{itemize}

Examples of values for this cost function are given in section \ref{s:WH} on an example of mean field control of water heaters.

\subsection{Definition of the problem}
We consider the mean field limit of the problem of controlling a large number of agents.
The objective is to minimize a cost function that depends on the mean-field interaction. Given a control policy $\alpha \in A$, the cost functional is defined as:
\begin{equation}
\label{e:J}
    J(\alpha) =   \int_0^T f(t, \mathbb{E}[p(t,X_t^{\alpha+\hat{\alpha}})]) dt + \mathbb{E}[G(\alpha,X^{\alpha+\hat{\alpha}}) ].
\end{equation}
{ We highlight that since we consider the mean field limit, the empirical distribution in \eqref{def_jn} in the coupling cost function $f$ is replaced above by the expectation of the individual consumption among the population of TCL.}
The optimal control problem is given by:
\begin{equation}
\label{e:minJ}
    \min_{\alpha \in A} \ J(\alpha),
\end{equation}
such that the dynamics of $X_t^{\alpha+\hat{\alpha}}$ follow the PDMP evolution described above. It is important to note here that the control parameter is not directly $I_t$ (the mode of heating for water heaters) but the part of the intensity that we can control $\alpha$.

\subsection{Assumptions}
We assume the following:

\textbf{(A1)} The continuous component is a.s. bounded: $\forall t,\theta_t\in[\theta_{0},\theta_{\infty}]$.

\textbf{(A2)} The function $p$ is in $C^1([0, T] \times \mathcal{X} )$ and $f : [0, T] \times \mathbb{R} \to \mathbb{R}$ is lower semi-continuous, strictly convex and twice differentiable with respect to the second variable with Lipschitz gradient uniformly in time.

\textbf{(A3)} $c \in C^1([0, T] \times \mathcal{X} )$, $g \in C^1(\mathcal{X})$ and $b$ is continuously differentiable with respect to the first variable.

\textbf{(A4)} The function $L: \mathbb{R} \to \mathbb{R}\cup \{+\infty\}$ is defined as:
\begin{equation*}
    L(x) :=
    \begin{cases}
        l(x), & \text{if } x > 0, \\
        0, & \text{if } x = 0, \\
        +\infty, & \text{otherwise},
    \end{cases}
\end{equation*}
where $l \in C^1(\mathbb{R}_+, \mathbb{R}_+)$ is an increasing strongly convex function bounded from above by a quadratic function.
We denote by \( H \) the convex conjugate of \( L \), i.e. \( H \) is defined for any \( x \in \mathbb{R} \) by
\begin{equation}
\label{e:Hconj}
H(x) := \sup_{y \in \mathbb{R}} \left( xy - L(y) \right).    
\end{equation}

These assumptions are verified for the use case defined in Section \ref{s:WH}. It is important to note that (A1) is weaker than the assumption made in \cite{seguret2024decentralized}, where the vector field $b$ must vanish on the edges of the interval. Assumption (A1) only requires having a zero density outside the domain $[\theta_0,\theta_\infty]$. The fact that the temperature of a TCL is bounded is a simple assumption to verify in practice: for water heaters, we can for example take $\theta_0$ to be the temperature of the water entering the heater, and $\theta_\infty$ to be the temperature reached if we were to heat the entire time interval from the maximum initial temperature. {  If the WH almost surely starts in this intervall, it will almost surely stay in this intervall}. These bounds are needed for theoretical reasons and should not be confused with the bounds we will impose in Section \ref{s:WH} (called $[\theta_{min},\theta_{max}]$), for reasons of quality of service or safety.

\subsection{Dualization}

First, we define a problem equivalent to the one defined in (\ref{e:minJ}).
With \( \alpha\in A \) and \( v \in L^2(0,T) \), we define:
\begin{equation*}
    \bar{J}(\alpha, v) =F(v) + \mathbb{E} \big[ G(\alpha, X^{\alpha+\hat{\alpha}}) \big]
\end{equation*}
where \( F \) is defined for any \( v \in L^2(0, T) \) by $F(v) := \int_0^T f(t, v(t)) dt$ and the equivalent problem is:
\begin{equation}
\label{e:minJv}
\begin{aligned}
        \min\limits_{\alpha\in A,v\in L^2(0,T)}&{\bar{J}(\alpha, v)} \\  \text{such that}& \ \mathbb{E} \big[ p(t, X^{\alpha+\hat{\alpha}}_t) \big] - v(t) = 0 \text{ a.e. on } [0, T]
\end{aligned}
\end{equation}
where the dynamics of $X_t^{\alpha+\hat{\alpha}}$ follow the PDMP evolution described above.

Let us introduce the Lagrangian 
\( \mathcal{L} : A \times L^2(0, T) \times L^2(0, T) \to \mathbb{R} \):

\begin{equation}
\begin{aligned}
        \mathcal{L}(\alpha, v, \lambda) :=& \mathbb{E} \big[ G(\alpha, X^{\alpha+\hat{\alpha}}) \big] + \int_0^T \mathbb{E} \big[ p(t, X^{\alpha+\hat{\alpha}}_t) \big] \lambda(t) dt\\&+ F(v) - \int_0^T v(t) \lambda(t) dt
\end{aligned}
\end{equation}
and the associated dual function 
\( W : L^2(0, T) \to \mathbb{R} \):
\begin{equation}
\label{e:WL1L2}
    W(\lambda) := \inf_{\alpha \in A,v \in L^2(0,T)} \mathcal{L}(\alpha,v, \lambda).
\end{equation}

The existence of a unique \( \bar{\lambda} \in L^2(0, T) \) such that
\[
    \bar{\lambda} = \arg\max_{\lambda \in L^2(0,T)} W(\lambda).
\]
is obtained in the Lemma 3.2 and Theorem 3.4 in \cite{seguret2024decentralized} (with a minor adaptation to take into account the additional Lipschitz continuous jump intensity $\hat{\alpha}$).
For a given $\lambda$, there exists a unique $ ( v[\lambda],\alpha[\lambda] )$ satisfying
    \begin{equation}
    \label{e:vmin}
        v[\lambda] = \arg\min_{v \in \mathcal{L}^2(0,T)} F(v) - \int_0^T v(t) \lambda(t).
    \end{equation}
\begin{equation}
\label{e:alpha}
\begin{aligned}
    \alpha_j& [\lambda](t, i,\theta) &\\&=\arg\min_{\alpha \in A}\mathbb{E} \big[ G(\alpha, X^{\alpha+\hat{\alpha}}) \big] + \int_0^T \mathbb{E} \big[ p(t, X^{\alpha+\hat{\alpha}}_t) \big] \lambda(t) dt\\&= H'\big( \varphi[\lambda](t, i,\theta) - \varphi[\lambda](t, j,\theta) \big)+ \hat{\alpha}_i\big(t,i,\theta)
    \end{aligned}
\end{equation}
where H is the convex conjugate of $L$ defined in (\ref{e:Hconj}) and \( \varphi[\lambda] \) is the unique function in \( C \) (defined in (1.1)), satisfying for any \( (t, i,\theta) \in [0, T] \times \mathcal{X} \),

\begin{equation}
\label{e:PhiInt}
\begin{aligned}
\varphi[\lambda](t, i,\theta) =& \int_t^T c\big( \tau, i, \phi_{t,\theta}(i, \tau) \big) + p\big( \tau, i, \phi_{t,\theta}(i, \tau) \big) \lambda(\tau)\\
&+\int_t^T \sum_{\substack{j \in I \\ j \neq i}} \Big[ -H \big( \varphi[\lambda](\tau, i, \phi_{t,\theta}(i, \tau)) \\ &-\varphi[\lambda](\tau, j, \phi_{t,\theta}(i, \tau)) \big) \Big]d\tau \\
&+ \int_t^T \sum_{\substack{j \in I \\ j \neq i}} 
\Big[ \hat{\alpha}_i\big(\tau,i,\phi_{t,\theta}(i,\tau)\big)\big( \varphi[\lambda](\tau, i, \phi_{t,\theta}(i, \tau)) 
\\&- \varphi[\lambda](\tau, j, \phi_{t,\theta}(i, \tau)) \big) \Big]d\tau + g\big( i, \phi_{t,\theta}(i, T) \big).
\end{aligned}
\end{equation} 

These results come from Theorem 2.1 and Proposition 4.1 in \cite{seguret2023}, with the addition of $\hat\alpha$.

Finally, adapting the result of Lemma 3.1 in \cite{seguret2024decentralized}, give us that \( W \) is Gâteaux differentiable in \( L^2(0, T) \). For any \( \lambda \in L^2(0, T) \), the Gâteaux differential of \( W \) at \( \lambda \) is denoted by \( D W(\lambda)(\cdot) \) and is given in the direction \( \mu \in L^2(0, T) \) by:
\begin{equation}\label{def_DW}
D W(\lambda)(\mu) = \int_0^T \left( E \left[ p(\tau, X^{\alpha[\lambda]+\hat{\alpha}}_\tau) \right] - v[\lambda](\tau) \right) \mu(\tau) \, d\tau    
\end{equation}

\section{ALGORITHM}
\label{s:Alg}
We therefore have in equation \ref{def_DW}, an expression for the differential that allows us to use a stochastic gradient descent algorithm by generating water heater trajectories to obtain an approximation of $E \left[ p(\tau, X^{\alpha[\lambda]+\hat{\alpha}}_\tau) \right]$. The step size of this algorithm is chosen as $\rho_k=\frac{a}{k+1}$, with $a$ chosen empirically to improve the speed of convergence.\vspace{-0.2cm}
\begin{algorithm}[H]
\caption{}
\begin{algorithmic}[1]
\STATE Initialize \( \lambda_0 \), set \( \{ \rho_k \} \) and \( M \in \mathbb{N}^* \)
\STATE \( k \gets 0 \)
\FOR{each \( k = 0, 1, \dots \)}
    \STATE \( v_k \gets v[\lambda_k] \) where \( v[\lambda_k] \) is defined in (\ref{e:vmin})
    \STATE \( \alpha_k \gets \alpha[\lambda_k] \) where \( \alpha[\lambda_k] \) is defined in (\ref{e:alpha})
    \STATE Generate \( M \) independent PDMP \( X_1^{\alpha_k+\hat{\alpha}}, \dots, X_M^{\alpha_k+\hat{\alpha}} \)
    \STATE \( U_{k+1} \gets \frac{1}{M} \sum_{j=1}^M p(\cdot, X_j^{ \alpha_k +\hat{\alpha}})-v_k \)
    \STATE \( \lambda_{k+1} \gets \lambda_k + \rho_k U_{k+1} \)
\ENDFOR
\end{algorithmic}
\end{algorithm}\vspace{-0.2cm}
Proofs of the convergence of such an algorithm can be found in \cite{seguret2024decentralized}.

{  To compute $\alpha_k$ at line 5 of this algorithm, it is necessary to discretize temperature and time and calculate $\varphi[\lambda_k]$ from (\ref{e:PhiInt}) via an equation backward in time.}
For this numerical scheme to be stable, the following necessary Courant–Friedrichs–Lewy condition \cite{courant1928uebe} must be met: $$\frac{Bdt}{d\theta}\leq 1 $$ where $B$ is the maximum value of the temperature derivative $b$ and $d\theta$ and $dt$ are the discretization sizes of temperature and time. Another necessary conditions for this numerical scheme to converge is: $dt\leq 1/(\| \alpha^\ast\|_\infty +\|\hat{ \alpha}\|_\infty )$ where $\alpha^\ast$ is the solution to the optimization problem. Although $\|\alpha^\ast\|_{\infty}$ is unknown, it is possible to find bounds (proof Lemma 3.9 of \cite{seguret2024decentralized}). These criteria enable us to know how to efficiently discretize space and time.

Computing this numerical scheme over the entire space $[\theta_{0},\theta_{\infty}] \times I$ can be time-consuming, particularly if the domain is large, especially due to the discretized continuous component. In the case of water heaters, where the aim is to prevent temperatures escaping from a desired temperature range $[\theta_{min},\theta_{max}]$, this scheme can only be calculated on this interval and its immediate proximity. Indeed, density outside this range is extremely low and has no impact on the resolution of the numerical scheme.

\section{NUMERICAL RESULTS}
\label{s:WH}

\subsection{Presentation of the use case}
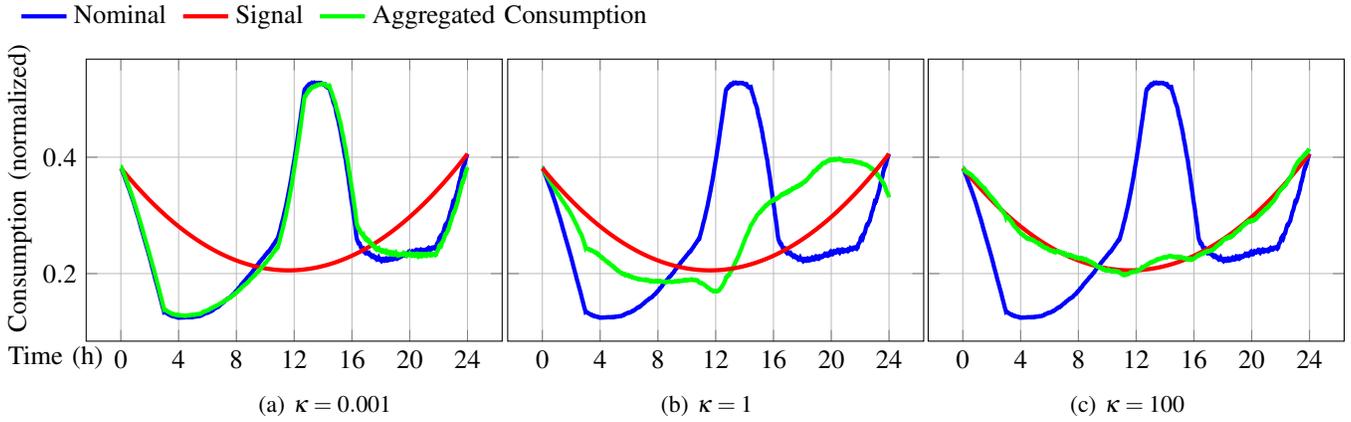
\begin{figure*}[t!]
\captionsetup[subfigure]{}
\hspace{-0.32cm}
\subcaptionbox{$\kappa=0.001$}{
\begin{tikzpicture}[scale=1]
\begin{axis}
[legend style={at={(0.6,1.23)}, anchor=north,draw=none},legend columns=5,
xtick  = {0,4,8,12,16,20,24},
xlabel={Time (h)},x label style ={at={(-0.21,0.017)},anchor=north west},
grid=major,
width=0.4\textwidth,height=0.3\textwidth,
ylabel={Consumption (normalized)},y label style ={at={(-0.21,0)},anchor=north west},legend entries={Nominal \quad ,Signal \quad ,Aggregated Consumption}]
\addplot [draw=blue,ultra thick] table[x index=0,y index=1]{DataFigures/rNom.txt};
\addplot [draw=red,ultra thick] table[x index=0,y index=1]{DataFigures/Signal.txt};
\addplot [draw=green,ultra thick, name path=f] table[x index=0,y index=1]{DataFigures/AggCons0.01.txt};
\end{axis}
\end{tikzpicture}}%
\hspace{-2.36cm}
\subcaptionbox{$\kappa=1$}{\begin{tikzpicture}[scale=1]
\begin{axis}
[grid=major,yticklabel=\empty,xtick  = {0,4,8,12,16,20,24},
width=0.4\textwidth,height=0.3\textwidth],
\addplot [draw=blue,ultra thick] table[x index=0,y index=1]{DataFigures/rNom.txt};
\addplot [draw=green,ultra thick, name path=f] table[x index=0,y index=1]{DataFigures/AggCons1.0.txt};
\addplot [draw=red,ultra thick] table[x index=0,y index=1]{DataFigures/Signal.txt};
\end{axis}
\end{tikzpicture}}%
\hspace{-0.28cm}
\subcaptionbox{$\kappa=100$}{\begin{tikzpicture}[scale=1]%
\begin{axis}
[grid=major,yticklabel=\empty,xtick  = {0,4,8,12,16,20,24},
width=0.4\textwidth,height=0.3\textwidth],
\addplot [draw=blue,ultra thick] table[x index=0,y index=1]{DataFigures/rNom.txt};
\addplot [draw=red,ultra thick] table[x index=0,y index=1]{DataFigures/Signal.txt};
\addplot [draw=green,ultra thick, name path=f] table[x index=0,y index=1]
{DataFigures/AggCons100.0.txt};
\end{axis}
\end{tikzpicture}}
\caption{Aggregated Consumptions for different values of $\kappa$ and $M=10^5$}
\vspace{-0.2cm}
\label{fig:3plots}
\end{figure*}

We consider a large number of water heaters that a central planner controls with a view to make their overall consumption approaches a certain signal $\{r_t\}_{t\in[0,T]}$, shown in red in Fig. \ref{fig:3plots}. We want to minimize the function $J$ defined in (\ref{e:J}) where:

\begin{itemize}
    \item The costs $c$ and $g$ are null and $\forall x\geq0,l(x)=\frac{x^2}{2}$.
    \item $p(t,X_t^{\alpha+\hat{\alpha}})$ is the power consumption associated to the PDMP $X_t^{\alpha+\hat{\alpha}}$ at time $t$ and \begin{equation}
    \label{e:Kappa}
        f(t,e)=\kappa(e-r_t)^2
    \end{equation} to penalize quadratically the distance to the signal. $\kappa> 0$ is a parameter used to balance the importance of signal tracking.
\end{itemize}

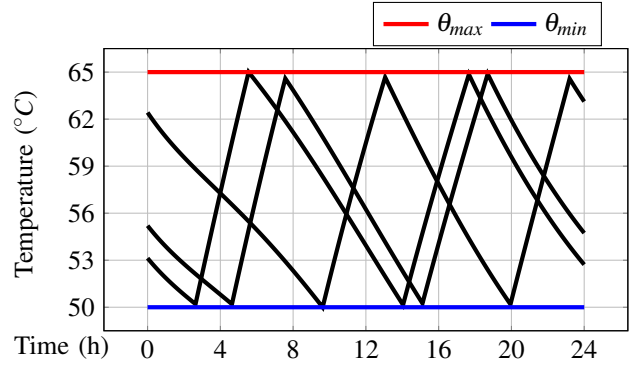
\begin{figure}[h]
\vspace{-0.2cm}
    \centering
    \begin{tikzpicture}[scale=0.43]
    \begin{axis}[xtick={-10,100},ytick={0,100},
    width=1\textwidth,height=0.58\textwidth,
    ylabel={Temperature ($^\circ C$)},y label style ={at={(-0.55,0.25)},anchor=north west}]
    \addplot [draw=black,ultra thick] table[x index=0,y index=1]{DataFigures/PDMPs.txt};
    \addplot [draw=black,ultra thick] table[x index=0,y index=2]{DataFigures/PDMPs.txt};
    \addplot [draw=black,ultra thick] table[x index=0,y index=3]{DataFigures/PDMPs.txt};
    \end{axis}

    \begin{axis}[legend style={at={(2.1,2.6)}, anchor=north east},legend columns=3,xtick  = {0,4,8,12,16,20,24},
    xlabel={Time (h)},x label style ={at={(-0.55,-0.08)},anchor=north west},ytick  = {50,53,56,59,62,65},
    grid=major,
    width=1\textwidth,height=0.58\textwidth,
    legend entries={$\theta_{max}$,$\theta_{min}$}]
    \addplot[red, ultra thick] coordinates {(0,65) (24,65)};
    \addplot[blue, ultra thick] coordinates {(0,50) (24,50)};
    \end{axis}
    \end{tikzpicture}
    \caption{Three trajectories of temperatures generated according to the nominal control}
    \vspace{-0.3cm}
    \label{fig:PDMPs}
\end{figure}

The time horizon is the day : $T=24h$ discretized every $dt=2$ minutes and the temperature set is discretized every $d\theta=1^\circ C$.

Water heaters are modeled by PMDPs where the temperatures evolve with the speed :
\begin{equation}
\frac{d\theta}{dt}=b(t,i,\theta)= \underbrace{i\sigma P}_{\textrm{heating effect}} - \underbrace{ \rho(\theta-\theta_{amb})}_{\textrm{loss effect}} - \underbrace{\epsilon (\theta-\theta_{in})}_{\textrm{drain effect}}
\label{eq:ODE}
\end{equation}
with $\rho$ the fraction of heat loss by minute, $\sigma$ the specific heat capacity of the volume of water, $P$, the heating power of the water heater, $\theta_{amb}$ the room temperature, $\theta_{in}$ the arriving water temperature and $\epsilon(t)$ the water drains at time $t$, which are considered similar for all water heaters. The time dependence of the $b$ velocity is due to the time-dependent drain term. The initial temperature $\theta_0$ is chosen with uniform probability between $\theta_{min}$ and $\theta_{max}$, and the initial mode is chosen ON (1) with probability 0.38 and OFF (0) with probability 0.62; these probability values are arbitrary.

For the water heater to operate correctly, temperatures must remain below a maximum temperature $\theta_{max}=65^\circ C$, and for the user's quality of service, temperatures must remain above a temperature $\theta_{min}=50^\circ C$ whenever possible. An additional intensity $\hat{\alpha}$, defined as follows, is therefore imposed to switch to heating mode when temperatures are too low and to switch off when they are too high: 

\begin{equation}
    \label{e:alphaH}
\begin{aligned}
        \hat{\alpha}_j(t,i,\theta) = \quad \ \quad \ \qquad \qquad\qquad \qquad\qquad \qquad\qquad \qquad  &\\ \left \{
\begin{array}{cc}
     12 & \text{if } \theta\geq\theta_{max}, i=1 \text{ and } j=0\\
     12\frac{\theta-\theta_{max}+d\theta}{d\theta} & \text{if } \theta\in[\theta_{max}-d\theta,\theta_{max}], i=1 \text{ and } j=0\\
     12 & \text{if } \theta\leq\theta_{min}, i=0 \text{ and } j=1\\
     12\frac{\theta_{min}+d\theta-\theta}{d\theta} & \text{if } \theta\in[\theta_{min},\theta_{min}+d\theta], i=1 \text{ and } j=0\\
     0 & \text{else}
\end{array}
\right.
\end{aligned}
\end{equation}
with $d\theta$, the discretization size of the temperature set. The value $12$ is chosen here empirically large enough to make PDMPs jump when arriving at the boundary while preserving the stability of the numerical scheme. It could be understood as follows: for a PDMP arriving at the boundary $\theta_{max}$ (or $\theta_{min}$), its next jump time follows the exponential distribution with parameter $12$. The expectancy of the jumping time when arriving at the boundary $\theta_{max}$ is thus $60/12=5$ minutes. Trajectories generated by this intensity, which we call \textit{nominal} in the following, are shown in Fig. \ref{fig:PDMPs}. These trajectories remain roughly within the predefined range.
\subsection{Signal Tracking}

Figure \ref{fig:3plots} shows overall consumption in blue for nominal control (without signal tracking). The signal in red, which we wish to approach, seeks to smooth consumption over the day. The green line shows aggregate consumption for several $ \kappa$ values ($\kappa$ is defined in (\ref{e:Kappa})). As expected, the higher $ \kappa$ goes, the closer the aggregated consumption is to the signal.

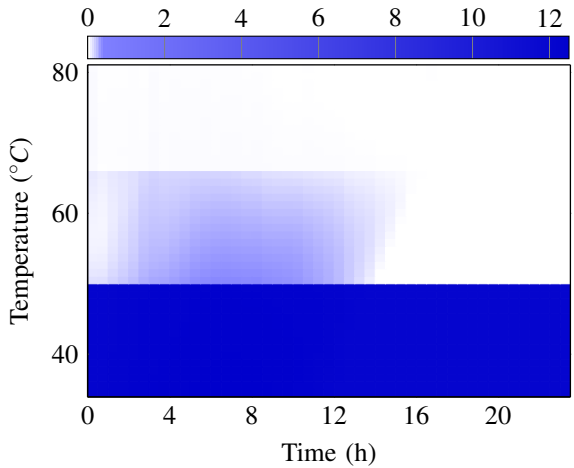
\begin{figure}[H]
    \centering
    \vspace{-0.5cm}
    \begin{tikzpicture}[scale=1]
\begin{axis}[view={0}{90},axis line style={line width=0.5pt, color=black},
             colormap={custom}{color(0)=(white) color(0.03)=(blue!50) color(1)=(blue!80!black)},
             colorbar horizontal,colorbar/width=3mm,
             colorbar style={anchor=north east, at={(6.4cm,5.1cm)},width=6.4cm,
             xticklabel pos=upper,scaled x ticks = false,
             xticklabel style={text width=4em,align=center,}},mesh/rows=45,
             mesh/cols=48,xtick={0,4,8,12,16,20,24},
             xlabel={Time (h)},
             ylabel={Temperature ($^\circ C$)},width=8 cm,height=6cm,tick style={major tick length=0pt}]
  \addplot3 [surf,shader=flat corner] table [row sep=newline] {DataFigures/alpha01K100.txt};
    \end{axis}
\end{tikzpicture}
    \caption{$\alpha+\hat{\alpha}$ for the jump from $0$ to $1$ (On to Off) for $\kappa=100$}
    \label{fig:alpha}
    \vspace{-0.4cm}
\end{figure}

If we look at the control $\alpha+\hat{\alpha}$ set up to get closer to the signal, we find that it is worth at least $12$ where water heaters are very likely to jump, so it is at least $\hat{\alpha}$. Typically, on Fig. \ref{fig:alpha}, we can see that for $\theta$ values below $\theta_{min}$, the value is at least $12$. Outside this range, lower values are found, but these are sufficient to approach the signal, with jumps favoured in the first half of the day to approach the signal, which is then higher than nominal consumption.

\subsection{Electricity cost}

\begin{figure}[H]
    \vspace{-0.2cm}
    \centering
    \begin{tikzpicture}[scale=0.42]
    \begin{axis}[legend style={draw=none,fill=none,at={(2.1,2.9)}, anchor=north east},legend columns=2,
xtick  = {0,4,8,12,16,20,24},
xlabel={Time (h)},x label style ={at={(-0.55,-0.08)},anchor=north west},
grid=major, ymin=-0.05,ymax=1.05,
width=0.9\textwidth,height=0.68\textwidth,
ylabel={Consumption (normalized)},y label style ={at={(-0.55,0.15)},anchor=north west},legend entries={Nominal, Aggregated Consumption}],
    \addplot [draw=blue,ultra thick] table[x index=0,y index=1]{DataFigures/rNom.txt};
    \addplot [draw=green,ultra thick] table[x index=0,y index=1]{DataFigures/rNomCost.txt};
    \end{axis}
    
    \begin{axis}[width=0.9\textwidth,height=0.68\textwidth,axis y line=right,axis x line=none,ylabel={Price (no unit)},y label style ={at={(2.7,1.8)},anchor=north west, rotate=180},legend entries={Electricity price},ymax=6,ymin=-6,legend style={draw=none,fill=none,at={(0.77,2.65)}, anchor=north east}]
    \addplot [draw=orange,ultra thick] table[x index=0,y index=1]{DataFigures/Cost.txt};
    \end{axis}
    \end{tikzpicture}
    \caption{An electricity cost is added}
    \label{fig:ElecCost}
    \vspace{-0.3cm}
\end{figure}
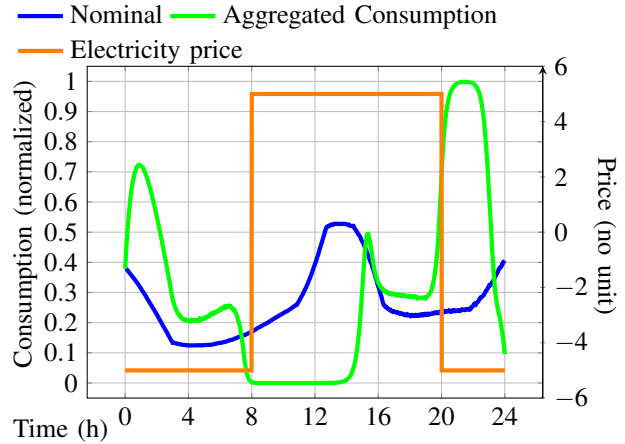

The term $c$ in (\ref{e:J}) corresponds to an individual cost for the agent. It is possible to use it to model an electricity price evolving throughout the day. On Fig. \ref{fig:ElecCost}, a given electricity cost (without signal tracking) is added and shown in orange. This cost of electricity would model a typical peak/off-peak situation where part of the day (here the middle of the day from 8am to 8pm) is penalized in terms of consumption price, while the rest of the day is not.

The corresponding aggregate consumption is large, as expected, before and after mid-day (when the price is low) and largely disappears in the middle of the day (when the price is high).  It does not disappear completely due to the need to stay within the temperature range (and therefore to heat when the temperature is too low). The algorithm chooses to heat water heaters to maximum before 8 a.m., restarting some at around 4 p.m. that would be too low in temperature, then at 8 p.m., most water heaters are very low in temperature and must be restarted. Consumption rises again at the end of the high-price window, as the algorithm penalizes changing modes too often. It therefore prefers to restart the heaters at the end of the window (and not have to switch them off afterwards) rather than have them heat up at the beginning of the window and have to switch them off afterwards, which incurs a cost due to the function $L$ introduced in (\ref{e:J}).

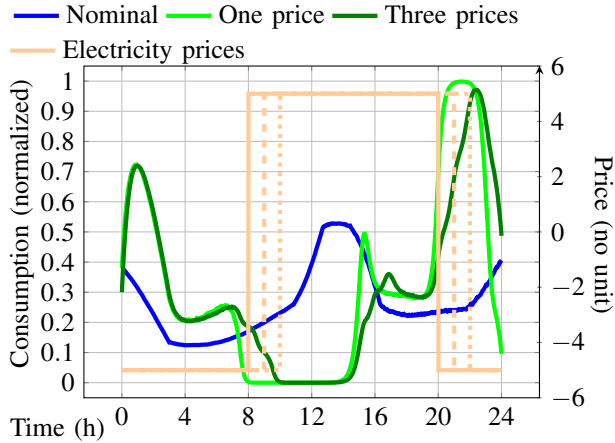
\begin{figure}[h!]
    \centering
    \begin{tikzpicture}[scale=0.42]
    \begin{axis}[legend style={draw=none,fill=none,at={(2.23,2.9)}, anchor=north east},legend columns=3,
xtick  = {0,4,8,12,16,20,24},
xlabel={Time (h)},x label style ={at={(-0.55,-0.08)},anchor=north west},
grid=major, ymin=-0.05,ymax=1.05,
width=0.9\textwidth,height=0.68\textwidth,
ylabel={Consumption (normalized)},y label style ={at={(-0.55,0.15)},anchor=north west},legend entries={Nominal, One price,Three prices}],
    \addplot [draw=blue,ultra thick] table[x index=0,y index=1]{DataFigures/rNom.txt};
    \addplot [draw=green,ultra thick] table[x index=0,y index=1]{DataFigures/rNomCost.txt};
    \addplot [draw=green!50!black,ultra thick] table[x index=0,y index=1]{DataFigures/rCost3.txt};
    \end{axis}
    \begin{axis}[width=0.9\textwidth,height=0.68\textwidth,axis y line=right,axis x line=none,ylabel={Price (no unit)},y label style ={at={(2.7,1.8)},anchor=north west, rotate=180},legend entries={Electricity prices},ymax=6,ymin=-6,legend style={draw=none,fill=none,at={(0.8,2.65)}, anchor=north east}]
    \addplot [draw=orange!40,ultra thick] table[x index=0,y index=1]{DataFigures/Cost.txt};
    \addplot[draw=orange!40,ultra thick,dashed] table[x index=0,y index=1]{DataFigures/Cost21.txt};
    \addplot[draw=orange!40,ultra thick,dotted] table[x index=0,y index=1]{DataFigures/Cost22.txt};
    \end{axis}
    \end{tikzpicture}
    \caption{Agents are divided into three classes, each with a different pricing structure, to smooth out the consumption}
    \label{fig:ElecCost2}
    \vspace{-0.2cm}
\end{figure}

In practice, we try to avoid sudden increases like the one at 8pm in Fig. \ref{fig:ElecCost}, for reasons of grid stability. For peak/off-peak tariffs, different classes of customers with different times of day are charged different prices. The reason for this can be seen in Fig. \ref{fig:ElecCost2}, where 3 different tariffs (high prices at 8 am/8 pm or 9 am/9 pm or 10 am/10 pm) are set up for 1/3 of the population each time. This helps to smooth out the consumption peak after 8pm. 

\section{Future works}


We feel that four possible paths could be explored: (i) integrating strong constraints that are often useful for modeling real-life situations (maximum power, smoothing slopes of aggregated consumption as in Fig. \ref{fig:ElecCost}, etc.). This could be done by integrating them into the Lagrangian directly, and would be easy to manage as the gradient according to the Lagrange multipliers would be expressed very similarly to those in this article (ii) Modify the differential equations governing PDMPs to integrate forced jumps, without the probabilistic aspect. Indeed, here the jump is made with a high probability, but an alternative approach could be to make it forced by modifying the continuity equation of the PDMP {  (iii) This framework requires knowledge of future information (such as water draw events, the reference signal, or electricity prices) to compute an optimal solution. A possible direction would therefore be to adapt it to an online setting, where predictions of these quantities are available at the beginning of the day. The algorithm would then be applied based on these forecasts, which would be progressively refined as new information is revealed at each time step, allowing for iterative updates and improved control performance (iv) The proposed approach relies on a key assumption: the number of agents must be sufficiently large. A sensitivity analysis with respect to the population size would offer valuable insights into the method's performance across different population scales.}

\bibliographystyle{plain}

\newpage

\end{document}